\documentclass[russcorr,a4paper,12pt]{iopart}
\usepackage{iopams}

\usepackage[russian]{babel}

\usepackage{graphicx}

\newtheorem{corollary}{Следствие}
\newtheorem{criterion}{Критерий}
\newtheorem{definition}{Определение}

\newtheorem{lemma}{Лемма}
\newtheorem{proposition}{Предложение}

\newenvironment{proof}[1][Доказательство]{\textbf{#1.} }{\ \rule{0.5em}{0.5em}}
\begin{document}

\title[Инварианты алгебр Ли размерности $n\leq 8$]{Инварианты коприсоединенного
представления алгебр Ли размерности $n\leq 8$\footnote{Работа
выполнена при финансовой поддержке UCM по проекту PR1-05/13283.}}
\author{Рутвиг Кампоамор-Штурсберг\footnote{Dpto. Geometr\'{\i}a y Topolog\'{\i}a, Fac. CC.
Matem\'aticas, Universidad Complutense de Madrid, Plaza de
Ciencias 3, E-28040 Mадрид. Эл. почта: rutwig@mat.ucm.es}}

\date{}

\begin{abstract}
Мы описываем инварианты корписоединненого представления для всех
вещественных алгебр Ли с нетривиальным разложением Леви до
размерности восемь.\newline УДК: 512.554.1, 539.12.01
\end{abstract}

\maketitle

\section{Введение. Постановка задачи}

Одной из важных характеристик алгебр Ли являются инвариантные
операторы коприсоединенного представления. Алгебраические
инварианты (классических) алгебр обычно называют инвариантами или
операторами Казимира. Эти инварианты играют важную роль в теории
представлений, а также во многих физических приложениях (массовые
формулы для мультиплетов, энергетические спектры, квантовые числа
и т.д.). Задача отыскания инвариантов Казимира (т.е.
полиномиальных инвариантов) для классических алгебр была решена в
ряде работ \cite{Ge,Pe,Pe1,Ra}. В частности, число независимых
инвариантов простой алгебры Ли равно рангу алгебры \cite{Ra}. В
этом случае, задача сводится к анализу центра универсальной
обертивающей алгебры $\mathcal{U}(\frak{g})$. В случае
неполупростых алгебр Ли инварианты коприсоединненого представления
могут быть рациональными или иррациональными. Имеются даже алгебры
Ли без инвариантов, как например $I\frak{gl}(n,\mathbb{R})$ или
двумерная неабелева алгебра $\frak{r}_{2}$. Это подразумевает, что
альтернативные методы должны быть развиты, чтобы вычуслуть
инварианты. Таким образом, вычислиение инвариантов требует
разработки иных методов.

В статьях \cite{De1,De2} были исследованы инварианты Казимира
алгебр Ли неоднородных групп и других алгебр,  имеющих важное
зчачение в физике. Инварианты трех-, четырех- и пятимерных алгебр
были вычислены в \cite{Wi}, а для разрешимых алгебр шестого
порядка в \cite{Nd} и \cite{C48}. Другие результаты об инвариантах
могут быть найдены в \cite{C35,C43,C52}.

\bigskip

Как известно, любая произвольная алгебра Ли $\frak{g}$ состоит из
двух частей: радикала $\frak{r}$ и факторалгебры $\frak{g}/
\frak{r}$. Фундаментальная теорема Леви-Мальцева дает полное
описание этого разложения (называемого разложением Леви). А
именно, любая алгебра Ли является полупростой суммой
$\frak{g}=\frak{s}\overrightarrow{\oplus}_{R}\frak{r}$ радикала
$\frak{r}$ и полупрямой подалгебры $\frak{s}$ (называется фактором
Леви), где $R$ --представление алгебры $\frak{s}$ т.ч.
\begin{equation}
R(X)*[Y,Z]=[R(X)*Y,Z]+[Y,R(X)*Z], \forall X\in\frak{s},\forall
Y,Z\in\frak{r}.
\end{equation}
В частности, для всех $X\in\frak{s}$ оператор $R$ -- внешнее
дифференцирование алгебры $\frak{r}$. \newline Это позволяет
свести задачу классификации алгебр Ли к следующим случаям:
\begin{enumerate}
\item Классификация всех полупростых алгебр Ли,

\item Классификация всех разрешимых алгебр Ли,

\item Классификация всех дифференцирований разрешимых алгебр в
сооветствии с классификацей простых алгебр Ли.
\end{enumerate}

Задача 1 решена полносью (см. например \cite{Vi}), но задачи 2 и 3
решены лишь частично (классификация вещественных разрешимых алгебр
малых размерностей известна, см. \cite{Mu,Tu}.) В работе
\cite{Tu3} была дана классификация неразложимых алгебр Ли типа
$\frak{g}=\frak{s}\overrightarrow{\oplus}_{R}\frak{r}$, которых
размерности не превышают девяти, а в работе \cite{Tu4}
классификаця алгебр Ли 10-го порядка специального вида.

В нашей работе найдены обобщенные инварианты коприсоединенного
представления неразложимых алгебр Ли $n$-го порядка с
нетривиальным разложением Леви, где $n\leq 8$.

\section{Аналитический метод и некоторые критерии}

Мы коротко напоминаем о стандартном методе, чтобы определять
инварианты коприсоединенного представления алгебры Ли $\frak{g}$,
в частности ее операторы Казимира. Пусть
$\left\{C_{ij}^{k}\right\}$ --- структурный тензор алгебры Ли
$\frak{g}$ в базисе $\left\{X_{1},..,X_{n}\right\}  $. Пусть
\begin{equation}
\widehat{X}_{i}=C_{ij}^{k}x_{k}\partial_{x_{j}}
\end{equation}
--- представление алгебры Ли дифференциальными операторами, где
$\left[  X_{i},X_{j}\right]  =C_{ij}^{k}X_{k}$ $\left(  1\leq
i<j\leq n,\;1\leq k\leq n\right)  $. Аналитическая функция
$F\left(x_{1},..,x_{n}\right)$ на $\frak{g}^{*}$ называется
инвариантом (коприсоединненого представления) алгебры $\frak{g}$
если $F(f)=F(ad^{*}_{g}(f))$ для любых $f\in \frak{g}$ и $g\in G$,
где G --- группа Ли, ассоциированная с алгеброй $\frak{g}$.

\begin{proposition}
Функция $F(x_{1},..,x_{n})$ на $\frak{g}^{*}$ является инвариантом
алгебры Ли тогда и только тогда, когда
\begin{equation}
\widehat{X}_{i}F(x_{1},..,x_{n})=C_{ij}^{k}x_{k}\frac{\partial
}{\partial x_{j}}F\left( x_{1},..,x_{n}\right) =0,\quad 1\leq
i\leq n.  \label{sys}
\end{equation}
\end{proposition}

Система (\ref{sys}) --- это система линейных однородных
дифференциальных уравнений в частных производных первого порядка.
Для алгебры Ли $\frak{g}$ существует полный набор из
$\mathcal{N}(\frak{g})$ инвариантов, где $\mathcal{N}(\frak{g})$
--- алгебра Ли определяется формулой
\begin{equation}
\mathcal{N}\left(  \frak{g}\right)  =\dim\,\frak{g}-{\rm rank}
\,\left(C_{ij}^{k}x_{k}\right).
\end{equation}

Для удобства изложения  мы вкратце перечисляем некоторые критерии,
используемые при решении задачи интегрирования системы (\ref{sys})
для алгебр Ли, которые мы рассматриваем в настоящей работе.

\begin{criterion}
Если $\frak{g}$ изоморфна производной алгебре
$\left[\frak{g},\frak{g}\right]$, то алгебра имеет только
полиномиальные инварианты.
\end{criterion}

Доказательство этой теоремы, как и рядa других, мы здесь не
рассматриваем. Его можно найти в \cite{AA}. Как можно доказать,
если (неполупростая) алгебра $\frak{g}$ удовлетворяет условию
$\frak{g}=\left[\frak{g},\frak{g}\right]$, то радикал алгебры
нильпотентен (см. например \cite{C35}). Для полупрямых сумм
специального вида имеет место следующая теорема \cite{C40}:

\begin{proposition}
Пусть $\frak{g}=\frak{s}\overrightarrow{\oplus}_{R}\frak{h}_{m}$
  $\left( m\geq1\right)$ ---алгебра Ли с одномерным центром, где $\frak{s=sl}\left(
2,\mathbb{R}\right)$ или $\frak{so}\left(3\right)$ и
$\frak{h}_{m}$ -- алгебра Гейзенберга размерности $(2m+1)$. Тогда
$\mathcal{N}\left(\frak{g}\right)=2$, и нецентральный инвариант
$C$ задается формулой
\begin{equation}
C=\sqrt{ \left| \left(
\begin{array}{ccccccc}
&  &  &  &  &  & x_{1} \\
&  &  &  &  &  & x_{2} \\
&  &  &  &  &  & x_{3} \\
&  & \left(C_{ij}^{k}x_{k}\right)  &  &  &  & \frac{1}{2}x_{4} \\
&  &  &  &  &  & : \\
&  &  &  &  &  & \frac{1}{2}x_{\dim{g}-1} \\
-x_{1} & -x_{2} & -x_{3} & -\frac{1}{2}x_{4} & ... &
-\frac{1}{2}x_{\dim{g}-1} & 0
\end{array}
\right) \right|}.
\end{equation}
\end{proposition}

Пусть $R$ ---представление алгебры Ли $\frak{s}$ и
$mult_{D_{0}}(R)$ -- кратность тривиального представления $D_{0}$
в $R$.

\begin{proposition}
Пусть $\frak{g}=\frak{s}\overrightarrow{\oplus}_{R}\frak{r}$. Если
$mult_{D_{0}}(R)=0$, то радикал нильпотентен. В частности,
$\mathcal{N}(\frak{g})>0$.
\end{proposition}

Доказательство можно найти в \cite{Tu3}. Тогда неразложимые
алгебры Ли, которые являются полупрямой суммой полупростой и
нильпотентной алгебр Ли удовлетворяют критерию 1
\footnote{Нильпотентные алгебры не удовлетворяют этому критерию,
но также имеют только инварианты Казимира \cite{Vi}. }. В
\cite{Tro} следующие результаты были получены, чтобы
проанализировать инварианты борелевских подалгебр полупростых
алгебр Ли:

\begin{definition}
Аналитическая функция $F(x_{1},..,x_{n})$ на $\frak{g}^{*}$
называется относи\-тельным инвариантом алгебры $\frak{g}$, если
\begin{equation}
\left[ F, X_{i}\right]=\lambda_{i}X_{i},
\end{equation}
для любых $X_{i}\in\frak{g}$, где $\lambda_{i}$ -- числовой
множитель.
\end{definition}

\begin{lemma}
Если $F$ и $G$ относительые инварианты алгебры $\frak{g}$ с
числовыми множителями $\lambda_{1}$ (сотв. $\lambda_{2}$), то
$F^{\lambda_{2}}G^{-\lambda_{1}}$ есть инвариант алгебры.
\end{lemma}

\begin{proof}
Пусть $X$ --- элемент алгебры $\frak{g}$. \
\begin{eqnarray}
\left[ X,F^{\lambda_{2}}G^{-\lambda_{1}}\right]=
XF^{\lambda_{2}}G^{-\lambda_{1}}-F^{\lambda_{2}}G^{-\lambda_{1}}X=
\left[X,F^{\lambda_{2}}\right]G^{-\lambda_{1}}-F^{\lambda_{2}}G^{\lambda_{1}}
\left[X,G^{\lambda_{1}}\right]G^{-\lambda_{1}} & \nonumber \\
=\lambda_{1}\lambda_{2}F^{\lambda_{2}}G^{-\lambda_{1}}-F^{\lambda_{2}}G^{-\lambda-{1}}(\lambda_{2}\lambda_{1}
G^{\lambda_{1}})G^{-\lambda_{1}}=0.&
\end{eqnarray}

\end{proof}

С помощью этой леммы можно изучать инварианты полупрямых сумм
$\frak{s}\overrightarrow{\oplus}_{R}\frak{r}$ полупростых и
разрешимых алгебр Ли.

\begin{proposition}
Пусть $\widehat{\frak{g}}$ --алгебра Ли и $\frak{g}$ -- подалгебра
коразмерности 1 такая, что имеет место следующее тождество:
\begin{equation}
\mathcal{N}(\widehat{\frak{g}})=\mathcal{N}(\frak{g})-1.
\end{equation}
Тогда всякий инвариант алгебры $\widehat{\frak{g}}$ является
инвариантом также и подалгебры $\frak{g}$. В частности, существует
элемент $X\notin \frak{g}$, такой что
\begin{equation}
\frac{\partial C}{\partial x}=0
\end{equation}
для каждого инварианта алгебры $\widehat{\frak{g}}$.
\end{proposition}

\begin{proof}
Шарп и др. показали в \cite{Sh}, что число необходимых внутренних
разметков индуцированного представления алгебры
$\widehat{\frak{g}}$  ассоциированное с редукцией
$\frak{g}\hookrightarrow \widehat{\frak{g}}$ равняется
\begin{equation}
n=\frac{1}{2}\left( \dim
\widehat{\frak{g}}-\mathcal{N}(\widehat{\frak{g}})-\dim \frak{g}-
\mathcal{N}(\frak{g})\right) +l^{\prime },  \label{ML}
\end{equation}
где $l^{\prime}$ -- число общих инвариантов, которие только
зависят от генераторов подалгебры $\frak{g}$. В этом случае
формула (\ref{ML}) дает:
\begin{equation}
n=\frac{1}{2}\left( -2\dim \mathcal{N}(\widehat{\frak{g}})\right)
+l^{\prime }\geq 0,
\end{equation}
тогда $l^{\prime}=\mathcal{N}(\widehat{\frak{g}})$. Из этого
следует, что всякий инвариант $F$ алгебры $\widehat{\frak{g}}$
зависит только от генераторов подалгебры $\frak{g}$. Это
подразумевает существование элемента $X\notin \frak{g}$, такой что
\begin{equation*}
\frac{\partial F}{\partial x}=0
\end{equation*}
для каждого инварианта алгебры $\widehat{\frak{g}}$.
\end{proof}

\begin{corollary}
Пусть $\frak{g}=\frak{s}\overrightarrow{\oplus}_{R}\frak{r}$ с
$mult_{D_{0}}(R)=1$. Пусть $X\in \frak{r}$ т.ч.
$\left[\frak{s},X\right]=0$. Если $\left[X,Z(\frak{n})\right]\neq
0$, где $\frak{n}$ -- нильрадикал алгебры $\frak{r}$, тогда
\begin{equation}
\frac{\partial F}{\partial x}=0
\end{equation}
для любого инварианта $F$ алгебры $\frak{g}$.
\end{corollary}

\begin{proposition}
Пусть $R$ -- неприводимое представление алгебры
$\frak{sl}(2,\mathbb{R})$ либо $\frak{so}(3)$. Если $\dim R>3$, то
существует максимальная система $F_{i}$ инвариантов Казимира
полупрямой суммы
$\frak{g}=\frak{s}\overrightarrow{\oplus}_{R}\dim{R}L_{1}$,
которые не зависят от переменных $x_{1},x_{2},x_{3}$,  отвечающих
фактору Леви.
\end{proposition}

Доказательства этих свойств могут быть найдены в \cite{C35,C48}.

\section{Иллюстративный пример}

Опишем, каким образом инварианты приводятся в таблицах 1-4.
Оставшиеся случаи исследуются аналогично. Пусть $L_{8,1}^{p}$
--- алгебра Ли, в которой скобка имеет вид:
\[
\begin{tabular}
[c]{llll}%
$\left[  X_{1},X_{2}\right]  =X_{3},$ & $\left[  X_{1},X_{3}\right]  =-X_{2},$%
& $\left[  X_{2},X_{3}\right]  =X_{1},$ & $\left[
X_{1},X_{5}\right]
=X_{6},$\\
$\left[  X_{1},X_{6}\right]  =-X_{5},$ & $\left[
X_{2},X_{4}\right]
=-X_{6},$ & $\left[  X_{2},X_{6}\right]  =X_{4},$ & $\left[  X_{3}%
,X_{4}\right]  =X_{5},$\\
$\left[  X_{3},X_{5}\right]  =-X_{4},$ & $\left[  X_{4},X_{8}\right]  =X_{4},$%
& $\left[  X_{5},X_{8}\right]  =X_{5},$ & $\left[
X_{6},X_{8}\right]
=X_{6},$\\
$\left[  X_{7},X_{8}\right]  =pX_{7},$ & $\left(  p\neq0\right) .$
&  &
\end{tabular}
\]
в базисе $\left\{X_{1},..,X_{8}\right\}$. В этом случае
$\mathcal{N}(L_{8,1})=2$. Система (\ref{sys}) примет вид системы 8
уравнений, из первых 7 уравнений которой следует, что инварианты
не зависят от $x_{8}$. Последнее уравнение этой системы имеет вид:
\begin{equation}
\widehat{X}_{8}F:=x_{4}\frac{\partial F}{\partial
x_{4}}+x_{5}\frac{\partial F}{\partial x_{5}}+x_{6}\frac{\partial
F}{\partial x_{6}}+px_{7}\frac{\partial F}{\partial x_{7}}=0.
\label{Eq}
\end{equation}
Алгебра $L_{8,1}^{p}$ содержит идеал 7-го порядка, который
изоморфен алгебре $L_{6,1}\oplus L_{1}$ (см. таблица 1). Эта
алгебра может иметь только полиномиальные инварианты. Фактически,
из разложимости следует, что
\begin{equation}
\mathcal{N}(L_{6,1}\oplus
L_{1})=\mathcal{N}(L_{6,1})+\mathcal{N}(L_{1}).
\end{equation}
Алгебра $L_{6,1}$ удовлетворяет критерию 1, следовательно, она
имеет только инварианты Казимира. Решая систему (\ref{sys})
подалгебры $L_{6,1}$ получаем инварианты:
\begin{eqnarray*}
I_{1}  & =x_{4}^{2}+x_{5}^{2}+x_{6}^{2}\\
I_{2}  & =x_{1}x_{4}+x_{2}x_{5}+x_{3}x_{6}\\
I_{3}  & =x_{7}.%
\end{eqnarray*}

Из предложения 4 следует, что $I_{1},I_{2},I_{3}$ являются
относительными инвариантами уравнения (\ref{Eq}). Числовые
множители имеют вид
\[
\widehat{X}_{8}\left(  I_{1}\right)
=-2I_{1},\;\widehat{X}_{8}\left( I_{2}\right)
=-I_{2},\;\widehat{X}_{8}\left(  I_{3}\right)  =-pI_{3},
\]
Из этого следует, что алгебра $L_{8,1}$ иммет рациональные
инварианты, которые согласно леммы 1 записываются в виде:
\begin{eqnarray*}
J_{1}  & =\frac{I_{1}^{p}}{I_{3}^{2}}=\frac{\left(  x_{4}^{2}+x_{5}^{2}%
+x_{6}^{2}\right)  ^{p}}{x_{7}^{2}}\\
J_{2}  & =\frac{I_{2}^{p}}{I_{3}}=\frac{\left(  x_{1}x_{4}+x_{2}x_{5}%
+x_{3}x_{6}\right)  ^{p}}{x_{7}}.
\end{eqnarray*}

\begin{table}
\caption{\label{Т1} Алгебры размерности $\leq 7$.}
\begin{tabular}
[c]{@{}lll}%
Алгебра & \textrm{ Структурные константы} & \textrm{Инварианты}
\\\hline
$L_{5,1}$ &
$C_{12}^{2}=2,C_{13}^{3}=-2,\;C_{23}^{1}=1,C_{14}^{4}=1,$ &
$I_{1}=x_{3}x_{4}^{2}-x_{1}x_{4}x_{5}-x_{2}x_{5}^{2}$\\
& $C_{25}^{4}=1,C_{34}^{5}=1,C_{15}^{5}=-1.$ & \\
$L_{6,1}$ &
$C_{23}^{1}=1,C_{12}^{3}=1,C_{13}^{2}=-1,C_{15}^{6}=1,$ &
$I_{1}=x_{4}^{2}+x_{5}^{2}+x_{6}^{2}$\\
& $C_{16}^{5}=-1,C_{24}^{6}=-1,C_{26}^{4}=1,C_{34}^{5}=1,$ & $I_{2}=x_{1}%
x_{4}+x_{2}x_{5}+x_{3}x_{6}$\\
& $C_{35}^{4}=-1.$ & \\
$L_{6,2}$ &
$C_{12}^{2}=2,C_{13}^{3}=-2,C_{23}^{1}=1,C_{14}^{4}=1,$ &
$I_{1}=2x_{2}x_{3}x_{6}+x_{5}^{2}x_{2}+x_{1}x_{4}x_{5}-x_{3}x_{4}^{2}+\frac
{1}{2}x_{1}^{2}x_{6}$\\
& $C_{25}^{4}=1,C_{34}^{5}=1,C_{15}^{5}=-1,C_{45}^{6}=1.$ & $I_{2}=x_{6}$\\
$L_{6,3}$ & $C_{12}^{2}=2,C_{13}^{3}=-2,C_{23}^{1}=1,C_{14}^{4}=1,$ & ---\\
& $C_{25}^{4}=1,C_{34}^{5}=1,C_{15}^{5}=-1,C_{j7}^{j}=1.$ & $(j=5,6)$\\
$L_{6,4}$ &
$C_{12}^{2}=2,C_{13}^{3}=-2,C_{23}^{1}=1,C_{14}^{4}=2,$ &
$I_{1}=x_{5}^{2}-4x_{4}x_{6}$\\
& $C_{16}^{6}=-2,C_{25}^{4}=2,C_{26}^{5}=1,C_{34}^{5}=1,$ & $I_{2}=x_{1}%
x_{5}+2x_{2}x_{6}-2x_{3}x_{4}$\\
& $C_{35}^{6}=2.$ & \\
$L_{7,1}$ &
$C_{23}^{1}=1,C_{12}^{3}=1,C_{13}^{2}=-1,\;C_{15}^{6}=1,$ &
$I_{1}=\left(  x_{4}^{2}+x_{5}^{2}+x_{6}^{2}\right)  \left(  x_{1}x_{4}%
+x_{2}x_{5}+x_{3}x_{6}\right)  ^{-2}$\\
& $C_{16}^{5}=-1,C_{24}^{6}=-1,C_{26}^{4}=1,C_{34}^{5}=1,$ & \\
& $C_{35}^{4}=-1,C_{j7}^{j}=1\;\left(  4\leq j\leq 6\right)  .$ & \\
$L_{7,2}$ & $C_{23}^{1}=1,C_{12}^{3}=1,C_{13}^{2}=-1,\;C_{14}^{7}=\frac{1}%
{2},$ & $I_{1}=x_{4}^{2}+x_{5}^{2}+x_{6}^{2}+x_{7}^{2}$\\
& $C_{15}^{6}=\frac{1}{2},C_{16}^{5}=-\frac{1}{2},C_{17}^{4}=-\frac{1}%
{2},C_{24}^{5}=\frac{1}{2},$ & \\
& $C_{25}^{4}=\frac{1}{2},C_{26}^{7}=\frac{1}{2},C_{27}^{6}=-\frac{1}%
{2},C_{34}^{6}=\frac{1}{2},$ & \\
&
$C_{35}^{7}=-\frac{1}{2},C_{36}^{4}=-\frac{1}{2},C_{37}^{5}=\frac{1}{2}.$
&
\\
$L_{7,3}$ &
$C_{12}^{2}=2,C_{13}^{3}=-2,C_{23}^{1}=1,C_{14}^{4}=1,$ &
$I_{1}=\left( x_{3}x_{4}^{2}-x_{1}x_{4}x_{5}-x_{2}x_{5}^{2}\right)
^{-p}x_{6}^{2}$\\
& $C_{15}^{5}=-1,C_{25}^{4}=1,C_{34}^{5}=1,C_{47}^{4}=1,$ & \\
& $C_{57}^{5}=1,C_{67}^{6}=p\;\left(  p\neq0\right)  .$ & \\
$L_{7,4}$ &
$C_{12}^{2}=2,C_{13}^{3}=-2,C_{23}^{1}=1,C_{14}^{4}=1,\;$ &
$I_{1}=(x_{5}^{2}x_{2}+x_{1}x_{4}x_{5}-x_{3}x_{4}^{2})x_{6}^{-1}+2x_{2}x_{3}+$\\
& $C_{15}^{5}=-1,C_{25}^{4}=1,C_{34}^{5}=1,C_{45}^{6}=1,$ & $+\frac{1}{2}x_{1}^{2}$\\
& $C_{47}^{4}=1,C_{57}^{5}=1,C_{67}^{6}=2.$ & \\
$L_{7,5}$ &
$C_{12}^{2}=2,C_{13}^{3}=-2,C_{23}^{1}=1,C_{14}^{4}=2,\;$ &
$I_{1}=\left(  x_{1}x_{5}+2x_{2}x_{6}-2x_{3}x_{4}\right)
^{2}\left(
x_{5}^{2}-4x_{4}x_{6}\right)  ^{-1}$\\
& $C_{16}^{6}=-2,C_{25}^{4}=2,C_{26}^{5}=1,C_{34}^{4}=1,$ & \\
& $C_{35}^{5}=2,C_{j7}^{j}=1\;\left(  j=,4,5,6\right)  .$ & \\
$L_{7,6}$ &
$C_{12}^{2}=2,C_{13}^{3}=-2,C_{23}^{1}=1,\;C_{14}^{4}=3,$ &
$I_{1}=27x_{4}^{2}x_{7}^{2}-18x_{4}x_{5}x_{6}x_{7}-x_{5}^{2}x_{6}^{2}+  $\\
& $C_{15}^{5}=1,C_{16}^{6}=-1,C_{17}^{7}=-3,C_{25}^{4}=3,\;$ &
$+4\left(x_{6}^{3}x_{4}+x_{7}x_{5}^{3}\right)$ \\
& $C_{26}^{5}=2,C_{27}^{6}=1,C_{34}^{5}=1,C_{35}^{6}=2,\;$ & \\
& $C_{36}^{7}=3.$ & \\
$L_{7,7}$ &
$C_{12}^{2}=2,C_{13}^{3}=-2,C_{23}^{1}=1,\;C_{14}^{4}=1,$ &
$I_{1}=x_{4}x_{7}-x_{5}x_{6}$\\
& $C_{15}^{5}=-1,\;C_{25}^{4}=1,C_{27}^{6}=1,C_{34}^{5}=1,$ & \\
& $C_{16}^{6}=1,C_{17}^{7}=-1,C_{36}^{7}=1.$ & \\
&  & \\\hline
\end{tabular}
\end{table}

\begin{table}
\caption{\label{Т3} Алгебры размерности $8$.}
\begin{tabular}
[c]{@{}lll}%
 Алгебра & \textrm{ Структурные константы} & \textrm{Инварианты}\\\hline
$L_{8,1}$ &
$C_{23}^{1}=1,C_{12}^{3}=1,C_{13}^{2}=-1,\;C_{15}^{6}=1,$ &
$I_{1}=\left(  x_{4}^{2}+x_{5}^{2}+x_{6}^{2}\right)  ^{p}x_{7}^{-2}$\\
& $C_{16}^{5}=-1,C_{24}^{6}=-1,C_{26}^{4}=1,C_{34}^{5}=1,$ &
$I_{2}=\left(
x_{1}x_{4}+x_{2}x_{5}+x_{3}x_{6}\right)  ^{p}x_{7}^{-1}$\\
& $C_{35}^{4}=-1,C_{j8}^{j}=1\;\left(  4\leq j\leq7\right)
,C_{78}^{7}=p.$ &
\\
$L_{8,2}$ & $C_{23}^{1}=1,C_{12}^{3}=1,C_{13}^{2}=-1,\;C_{14}^{7}=\frac{1}%
{2},$ & $I_{1}=16(x_{4}x_{5}+x_{6}x_{7})+\left(  x_{4}^{2}+x_{5}^{2}+x_{6}^{2}+x_{7}^{2}\right)  $\\
& $C_{15}^{6}=\frac{1}{2},C_{16}^{5}=-\frac{1}{2},C_{17}^{4}=-\frac{1}%
{2},C_{24}^{5}=\frac{1}{2},$ & $+16(x_{1}^{2}+x_{2}^{2}+x_{3}^{2})x_{8}^{2}+16x_{2}x_{8}(x_{5}x_{6}-x_{4}x_{7}$\\
&
$C_{25}^{4}=\frac{1}{2},C_{26}^{7}=\frac{1}{2},C_{27}^{6}=-\frac{1}
{2},C_{34}^{6}=\frac{1}{2},$ & $+8x_{3}x_{8}(x_{4}^{2}-x_{5}^{2}+x_{6}^{2}-x_{7}^{2})$\\
& $C_{35}^{7}=-\frac{1}{2},C_{36}^{4}=-\frac{1}{2},C_{37}^{5}=\frac{1}%
{2},C_{45}^{8}=1,$ & \\
& $C_{67}^{8}=-1.$ & $I_{2}=x_{8}$\\
$L_{8,3}$ & $C_{23}^{1}=1,C_{12}^{3}=1,C_{13}^{2}=-1,\;C_{14}^{7}=\frac{1}%
{2},$ & \\
& $C_{15}^{6}=\frac{1}{2},C_{16}^{5}=-\frac{1}{2},C_{17}^{4}=-\frac{1}%
{2},C_{24}^{5}=\frac{1}{2},$ & \\
& $C_{25}^{4}=\frac{1}{2},C_{26}^{7}=\frac{1}{2},C_{27}^{6}=-\frac{1}%
{2},C_{34}^{6}=\frac{1}{2},$ & \\
& $C_{35}^{7}=-\frac{1}{2},C_{36}^{4}=-\frac{1}{2},C_{37}^{5}=\frac{1}%
{2},C_{48}^{4}=1,$ & \\
& $C_{58}^{5}=1,C_{68}^{6}=1,C_{78}^{7}=1.$ & \\
$L_{8.4}^{p}$ &
$C_{23}^{1}=1,C_{12}^{3}=1,C_{13}^{2}=-1,\;C_{14}^{7}=\frac
{1}{2},$ & $I_{1}=x_{4}^{2}+x_{5}^{2}+x_{6}^{2}+x_{7}^{2}$\\
& $C_{15}^{6}=\frac{1}{2},C_{16}^{5}=-\frac{1}{2},C_{17}^{4}=-\frac{1}%
{2},C_{24}^{5}=\frac{1}{2},$ & $I_{2}=\left(
x_{4}^{2}+x_{6}^{2}\right) \left(  x_{8}-2x_{3}\right)  +\left(
x_{5}^{2}+x_{7}^{2}\right)x_{8}+ $\\
& $C_{25}^{4}=\frac{1}{2},C_{26}^{7}=\frac{1}{2},C_{27}^{6}=-\frac{1}%
{2},C_{34}^{6}=\frac{1}{2},$ & $+4\left(  x_{2}x_{4}x_{7}-x_{1}x_{4}%
x_{5}-x_{2}x_{5}x_{6}-x_{1}x_{6}x_{7}\right) +$\\
& $C_{35}^{7}=-\frac{1}{2},C_{36}^{4}=-\frac{1}{2},C_{37}^{5}=\frac{1}%
{2},C_{48}^{48}=p,$ & $+2x_{3}\left( x_{5}^{2}+x_{7}^{2}\right)$\\
& $C_{58}^{5}=p,C_{68}^{6}=p,C_{78}^{7}=p,C_{48}^{6}=-1,$ & \\
& $C_{58}^{7}=-1,C_{68}^{4}=1,C_{78}^{5}=1.$ & \\
$L_{8,5}$ & $C_{23}^{1}=1,C_{12}^{3}=1,C_{13}^{2}=-1,C_{14}^{7}=\frac{1}{2},$%
& $I_{1}=x_{8}\left(
x_{7}^{2}+x_{6}^{2}-x_{6}^{2}-8x_{4}^{2}\right)
+6\left(  x_{6}^{2}-x_{7}^{2}\right)  +$\\
&
$C_{15}^{6}=-\frac{1}{2},C_{16}^{5}=2,C_{16}^{8}=-1,C_{17}^{4}=-2,$
&
$+\frac{2}{9}x_{8}^{3}-12x_{4}x_{6}x_{7}$\\
& $C_{18}^{6}=3,C_{24}^{6}=\frac{1}{2},C_{25}^{7}=\frac{1}{2},C_{26}^{4}=-2,$%
& $I_{2}=12\left(  x_{4}^{2}+x_{5}^{2}\right)  +3\left(  x_{6}^{2}+x_{7}%
^{2}\right)  +x_{8}^{2}$\\
& $C_{27}^{5}=-2,C_{27}^{8}=-1,C_{28}^{7}=3,C_{34}^{5}=2,$ & \\
& $C_{35}^{4}=-2,C_{36}^{7}=1,C_{37}^{6}=-1.$ & \\
$L_{8,6}$ &
$C_{12}^{2}=2,C_{13}^{3}=-2,C_{23}^{1}=1,C_{14}^{4}=1,$ &
$I_{1}=x_{8}$\\
& $C_{15}^{5}=-1,C_{25}^{4}=1,C_{34}^{5}=1,C_{67}^{8}=1.$ & $I_{2}=2x_{1}%
x_{4}x_{5}-2x_{3}x_{4}^{2}+2x_{2}x_{5}^{2}+4x_{2}x_{3}x_{8}+$\\
 & & $+x_{1}^{2}x_{8}$\\
$L_{8,7}^{p,q}$ &
$C_{12}^{2}=2,C_{13}^{3}=-2,C_{23}^{1}=1,C_{14}^{4}=1,$ &
$I_{1}=\left(  x_{3}x_{4}^{2}-x_{1}x_{4}x_{5}-x_{2}x_{5}^{2}\right)  ^{p}%
x_{6}^{-2}$\\
pq$\neq0$ &
$C_{15}^{5}=-1,C_{25}^{4}=1,C_{34}^{5}=1,C_{48}^{4}=1,$ &
$I_{2}=\left(  x_{3}x_{4}^{2}-x_{1}x_{4}x_{5}-x_{2}x_{5}^{2}\right)  ^{q}%
x_{7}^{-2}$\\
& $C_{58}^{5}=1,C_{68}^{6}=p,C_{78}^{7}=q.$ & \\
$L_{8,8}^{p}$ &
$C_{12}^{2}=2,C_{13}^{3}=-2,C_{23}^{1}=1,C_{14}^{4}=1,$ &
$I_{1}=\left(  x_{3}x_{4}^{2}-x_{1}x_{4}x_{5}-x_{2}x_{5}^{2}\right)  ^{p}%
x_{6}^{-2}$\\
p$\neq0$ & $C_{15}^{5}=-1,C_{25}^{4}=1,C_{34}^{5}=1,C_{48}^{4}=1,$
&
$I_{2}=\left(  p\ln x_{7}-x_{6}\ln x_{6}\right)  \left(  px_{6}\right)  ^{-1}%
$\\
& $C_{58}^{5}=1,C_{68}^{6}=p,C_{78}^{6}=1,C_{78}^{7}=p.$ & \\
$L_{8.8}^{0}$ &
$C_{12}^{2}=2,C_{13}^{3}=-2,C_{23}^{1}=1,C_{14}^{4}=1,$ &
$I_{1}=x_{6}$\\
& $C_{15}^{5}=-1,C_{25}^{4}=1,C_{34}^{5}=1,C_{48}^{4}=1,$ & $I_{2}%
=2x_{7}-x_{6}\ln\left(
x_{3}x_{4}^{2}-x_{1}x_{4}x_{5}-x_{2}x_{5}^{2}\right)
$\\
& $C_{58}^{5}=1,C_{78}^{6}=1.$ & \\
&  & \\\hline
\end{tabular}
\end{table}

\begin{table}
\caption{\label{Т4} Алгебры размерности $8$. (продолжение)}
\begin{tabular}
[c]{@{}lll}%
 Алгебра & \textrm{ Структурные константы} & \textrm{Инварианты} \\\hline
$L_{8,9}^{p,q}$ &
$C_{12}^{2}=2,C_{13}^{3}=-2,C_{23}^{1}=1,C_{14}^{4}=1,$ &
$I_{1}=\left(  x_{7}-ix_{6}\right)  ^{p-iq}\left(
x_{7}+ix_{6}\right)
^{p+iq}$\\
$q\neq0$ & $C_{15}^{5}=-1,C_{25}^{4}=1,C_{34}^{5}=1,C_{48}^{4}=1,$
& $I_{2}=\left(x_{3}x_{4}^{2}-x_{1}x_{4}x_{5}-x_{5}^{2}x_{2}\right)  ^{p^{2}+q^{2}}\times$\\
& $C_{58}^{5}=1,C_{68}^{6}=p,C_{68}^{7}=-q,C_{78}^{7}=q,$ &
$\left(  x_{7}-ix_{6}\right)  ^{2\left(  iq-p\right)
}$ \\
& $C_{78}^{7}=p.$ & \\
$L_{8,10}^{p}$ &
$C_{12}^{2}=2,C_{13}^{3}=-2,C_{23}^{1}=1,C_{14}^{4}=1,$ &
$I_{1}=\left(  2x_{2}x_{3}x_{6}+x_{5}^{2}x_{2}+x_{1}x_{4}x_{5}-x_{3}x_{4}%
^{2}\right)  x_{6}^{-1}$\\
& $C_{15}^{5}=-1,C_{25}^{4}=1,C_{34}^{5}=1,C_{48}^{4}=1,$ & $\frac{1}{2}x_{1}^{2}$\\
& $C_{58}^{5}=1,C_{68}^{6}=2,C_{78}^{7}=p.$ & $I_{2}=x_{7}%
^{2}x_{6}^{-p}$ \\
$L_{8,11}$ &
$C_{12}^{2}=2,C_{13}^{3}=-2,C_{23}^{1}=1,C_{14}^{4}=1,$ &
$I_{1}=\left(  2x_{2}x_{3}x_{6}+x_{5}^{2}x_{2}+x_{1}x_{4}x_{5}-x_{3}x_{4}%
^{2}\right)  x_{6}^{-1}$\\
& $C_{15}^{5}=-1,C_{25}^{4}=1,C_{34}^{5}=1,C_{48}^{4}=1,$ &
$+\frac{1}{2}x_{1}^{2}x_{6}$\\
& $C_{58}^{5}=1,C_{68}^{6}=2,C_{78}^{6}=1,C_{78}^{7}=p.$ &
$I_{2}=\left(
2x_{7}-x_{6}\ln x_{6}\right)  x_{6}^{-1}$ \\
$L_{8,12}$ &
$C_{12}^{2}=2,C_{13}^{3}=-2,C_{23}^{1}=1,C_{14}^{4}=2,$ &
$I_{1}=\left(  x_{5}^{2}-4x_{4}x_{6}\right)  ^{p}x_{7}^{-2}$\\
& $C_{16}^{6}=-2,C_{25}^{4}=2,C_{26}^{5}=1,C_{34}^{4}=1,$ &
$I_{2}=\left(
x_{1}x_{5}-2x_{3}x_{4}+2x_{2}x_{6}\right)  ^{p}x_{7}^{-2}$\\
& $C_{35}^{5}=2,C_{48}^{4}=1,C_{58}^{5}=1,C_{68}^{6}=1,$ & \\
& $C_{78}^{7}=p.$ & \\
$L_{8,13}^{\varepsilon}$ & $C_{12}^{2}=2,C_{13}^{3}=-2,C_{23}^{1}=1,C_{14}%
^{4}=1,$ & $I_{1}=x_{8}$\\
& $C_{15}^{5}=-1,C_{25}^{4}=1,C_{34}^{5}=1,C_{16}^{6}=1,$ &
$I_{2}=\varepsilon
x_{8}^{2}\left(  4x_{2}x_{3}+x_{1}^{2}\right)  +2x_{2}x_{8}\left(  x_{7}%
^{2}+\varepsilon x_{5}^{2}\right) $ \\
& $C_{17}^{7}=-1,C_{27}^{6}=1,C_{36}^{7}=1,C_{45}^{8}=1,$ & $-2x_{3}%
x_{8}\left(  \varepsilon x_{4}^{2}+x_{6}^{2}\right)
+2x_{6}x_{7}\left(
x_{1}x_{8}+x_{4}x_{5}\right)  $\\
& $C_{67}^{8}=\varepsilon.$ & $-x_{4}^{2}x_{7}^{2}$ \\
$L_{8,14}$ &
$C_{12}^{2}=2,C_{13}^{3}=-2,C_{23}^{1}=1,C_{14}^{4}=1,$ &
$I_{1}=x_{4}x_{7}-x_{5}x_{6}$\\
& $C_{15}^{5}=-1,C_{25}^{4}=1,C_{34}^{5}=1,C_{16}^{6}=1,$ & $I_{2}=x_{1}%
x_{4}x_{5}+x_{5}x_{6}x_{8}-x_{4}x_{7}x_{8}+x_{2}x_{5}^{2}$\\
& $C_{17}^{7}=-1,C_{27}^{6}=1,C_{36}^{7}=1,C_{48}^{4}=1,$ & \\
& $C_{58}^{5}=1.$ & $-x_{3}x_{4}^{2}$ \\
$L_{8,15}$ &
$C_{12}^{2}=2,C_{13}^{3}=-2,C_{23}^{1}=1,C_{14}^{4}=1,$ &
$I_{1}=x_{4}x_{7}-x_{5}x_{6}-\frac{1}{2}x_{8}^{2}$\\
& $C_{15}^{5}=-1,C_{25}^{4}=1,C_{34}^{5}=1,C_{16}^{6}=1,$ & $I_{2}=x_{1}%
x_{4}x_{5}+x_{5}x_{6}x_{8}-x_{4}x_{7}x_{8}+x_{2}x_{5}^{2}$\\
& $C_{17}^{7}=-1,C_{27}^{6}=1,C_{36}^{7}=1,C_{67}^{8}=1,$ & $-\frac{1}{3}x_{8}^{3}-x_{3}x_{4}^{2}$\\
& $C_{68}^{4}=1,C_{78}^{5}=1.$ & \\
$L_{8,16}$ & $C_{12}^{2}=2,C_{13}^{3}=-2,C_{23}^{1}=1,C_{14}^{4}=1,$ & \\
& $C_{15}^{5}=-1,C_{25}^{4}=1,C_{34}^{5}=1,C_{16}^{6}=1,$ & \\
& $C_{17}^{7}=-1,C_{27}^{6}=1,C_{36}^{7}=1,C_{48}^{4}=1,$ & \\
& $C_{58}^{5}=1,C_{68}^{4}=1,C_{68}^{6}=1,C_{78}^{5}=1,$ & \\
& $C_{78}^{7}=1.$ & \\
$L_{8,17}^{p}$ & $C_{12}^{2}=2,C_{13}^{3}=-2,C_{23}^{1}=1,C_{14}^{4}=1,$ & \\
$p\neq-1$ & $C_{15}^{5}=-1,C_{25}^{4}=1,C_{34}^{5}=1,C_{16}^{6}=1,$ & \\
& $C_{17}^{7}=-1,C_{27}^{6}=1,C_{36}^{7}=1,C_{48}^{4}=1,$ & \\
& $C_{58}^{5}=1,C_{68}^{6}=p,C_{78}^{7}=p.$ & \\
$L_{8,17}^{-1}$ &
$C_{12}^{2}=2,C_{13}^{3}=-2,C_{23}^{1}=1,C_{14}^{4}=1,$ &
$I_{1}=x_{4}x_{7}-x_{5}x_{6}$\\
& $C_{15}^{5}=-1,C_{25}^{4}=1,C_{34}^{5}=1,C_{16}^{6}=1,$ & $I_{2}%
=x_{1}\left(  x_{4}x_{7}+x_{5}x_{6}\right)  +2x_{2}x_{5}x_{7}-2x_{3}x_{4}%
x_{6}+$\\
& $C_{17}^{7}=-1,C_{27}^{6}=1,C_{36}^{7}=1,C_{48}^{4}=1,$ &
$x_{8}\left(
x_{4}x_{7}-x_{5}x_{6}\right)  $\\
& $C_{58}^{5}=1,C_{68}^{6}=-1,C_{78}^{7}=-1.$ & \\\hline
\end{tabular}
\end{table}

\begin{table}
\caption{\label{Т5} Алгебры размерности $8$. (продолжение)}
\begin{tabular}
[c]{@{}lll}%
Алгебра & \textrm{ Структурные константы} & \textrm{Инварианты}
\\\hline $L_{8,18}^{p}$ &
$C_{12}^{2}=2,C_{13}^{3}=-2,C_{23}^{1}=1,C_{14}^{4}=1,$ &
\\
$p\neq0$ & $C_{15}^{5}=-1,C_{25}^{4}=1,C_{34}^{5}=1,C_{16}^{6}=1,$ & \\
& $C_{17}^{7}=-1,C_{27}^{6}=1,C_{36}^{7}=1,C_{48}^{4}=p,$ &
$x_{8}\left(
x_{4}x_{7}-x_{5}x_{6}\right)  $\\
& $C_{48}^{6}=-1,C_{58}^{5}=p,C_{58}^{7}=-1,C_{68}^{4}=1,$ & \\
& $C_{68}^{6}=p,C_{78}^{5}=1,C_{78}^{7}=p$ & \\
$L_{8,18}^{0}$ &
$C_{12}^{2}=2,C_{13}^{3}=-2,C_{23}^{1}=1,C_{14}^{4}=1,$ &
$I_{1}=x_{4}x_{7}-x_{5}x_{6}$\\
& $C_{15}^{5}=-1,C_{25}^{4}=1,C_{34}^{5}=1,C_{16}^{6}=1,$ & $I_{2}%
=x_{1}\left(  x_{4}x_{7}+x_{5}x_{6}\right)  +x_{2}\left(  x_{5}^{2}+x_{7}%
^{2}\right)  -x_{3}x_{4}^{2}+$\\
& $C_{17}^{7}=-1,C_{27}^{6}=1,C_{36}^{7}=1,C_{48}^{6}=-1,$ &
$+x_{8}\left(
x_{5}x_{6}-x_{4}x_{7}\right) -x_{3}x_{6}^{2} $\\
& $C_{58}^{7}=-1,C_{68}^{4}=1,C_{78}^{5}=1.$ & \\
$L_{8,19}$ &
$C_{12}^{2}=2,C_{13}^{3}=-2,C_{23}^{1}=1,C_{14}^{4}=3,$ &
$I_{1}=x_{8}$\\
& $C_{15}^{5}=1,C_{16}^{6}=-1,C_{17}^{7}=-3,C_{25}^{4}=3,$ &
$I_{2}=12\left(
x_{2}x_{3}x_{8}^{2}+x_{2}x_{5}x_{7}x_{8}-x_{3}x_{4}x_{6}x_{8}\right)
+$\\
& $C_{26}^{5}=2,C_{27}^{6}=1,C_{34}^{5}=1,C_{35}^{6}=2,$ &
$+4\left(
x_{4}x_{6}^{3}+x_{3}x_{5}^{2}x_{8}+x_{5}^{3}x_{7}-x_{2}x_{6}^{2}x_{8}\right)
+$\\
& $C_{36}^{7}=3,C_{47}^{8}=1,C_{56}^{8}=-3.$ & $+18\left(  x_{1}x_{4}%
x_{7}x_{8}-x_{4}x_{5}x_{6}x_{7}\right)  +27x_{4}^{2}x_{7}^{2}+$\\
 & & $+3x_{1}^{2}x_{8}^{2}-2x_{1}x_{5}x_{6}x_{8}-x_{5}^{2}%
x_{6}^{2}$\\
$L_{8,20}$ & $C_{12}^{2}=2,C_{13}^{3}=-2,C_{23}^{1}=1,C_{14}^{4}=3,$ & \\
& $C_{15}^{5}=1,C_{16}^{6}=-1,C_{17}^{7}=-3,C_{25}^{4}=3,$ & \\
& $C_{26}^{5}=2,C_{27}^{6}=1,C_{34}^{5}=1,C_{35}^{6}=2,$ & \\
& $C_{36}^{7}=3,C_{48}^{4}=1,C_{58}^{5}=1,C_{68}^{6}=1,$ & \\
& $C_{78}^{7}=1.$ & \\
$L_{8,21}$ &
$C_{12}^{2}=2,C_{13}^{3}=-2,C_{23}^{1}=1,C_{14}^{4}=4,$ &
$I_{1}=3x_{5}x_{7}-x_{6}^{2}-12x_{4}x_{8}$\\
& $C_{15}^{5}=2,C_{17}^{7}=-2,C_{18}^{8}=-4,C_{25}^{4}=4,$ & $I_{2}%
=-2x_{6}^{3}+9x_{5}x_{6}x_{7}+72x_{4}x_{6}x_{8}-27x_{4}x_{7}^{2}$\\
& $C_{26}^{5}=3,C_{27}^{6}=2,C_{28}^{7}=1,C_{34}^{5}=1,$ & $-27x_{5}^{2}x_{8}$\\
& $C_{35}^{6}=2,C_{36}^{7}=3,C_{37}^{8}=4.$ & \\
$L_{8,22}$ &
$C_{12}^{2}=2,C_{13}^{3}=-2,C_{23}^{1}=1,C_{14}^{4}=2,$ &
$I_{1}=x_{5}^{2}-4x_{4}x_{6}$\\
& $C_{16}^{6}=-2,C_{17}^{7}=1,C_{18}^{8}=-1,C_{25}^{4}=2,$ & $I_{2}=x_{4}%
x_{8}^{2}+x_{6}x_{7}^{2}-x_{5}x_{7}x_{8}$\\
& $C_{26}^{5}=1,C_{28}^{7}=1,C_{34}^{5}=1,C_{35}^{6}=2,$ & \\
& $C_{37}^{8}=1.$ & \\\hline
\end{tabular}
\end{table}

\newpage

\section*{Список литературы }

\end{document}